\documentclass[preprint, 3p, authoryear]{elsarticle}
\usepackage{graphicx}
\usepackage{amssymb,amsthm,amsmath}
\usepackage{xcolor,paralist,hyperref,fancyhdr,etoolbox}
\usepackage{subfig}
\usepackage{mathtools}
\usepackage{cite}
\usepackage{ulem}
\usepackage{enumitem}
\usepackage{soul}
\usepackage{cleveref}
\usepackage{float}
\journal{}

\begin{document}

\begin{frontmatter}
\title{Markovian Embedding of Nonlinear Memory via Spectral Representation}
%\title{\rv{Markovian embedding of distributed delay systems with nonlinear memory effects using spectral representation}}%Markovian embedding of some nonlocal equations using spectral representation}
\author[inst1,inst2]{Divya Jaganathan}
\ead{djaganathan@seas.harvard.edu}
\affiliation[inst1]{organization={International Centre for Theoretical Sciences, Tata Institute of Fundamental Research},
            city={Bengaluru},
            postcode={560089}, 
            country={India}}
\affiliation[inst2]{organization={School of Engineering and Applied Sciences, Harvard University},city={Cambridge}, postcode={MA 02138}, country={USA}}
\affiliation[inst3]{organization={School of Computer and Mathematical Sciences, University of Adelaide},
            postcode={5005},
            country={Australia}}
\affiliation[inst4]{organization={Rudolf Peierls Centre for Theoretical Physics, Parks Road,
University of Oxford},
            city={Oxford},
            postcode={OX1 3PU},
            country={United Kingdom}}    

\author[inst3,inst4]{Rahil N. Valani}

\begin{abstract}
Differential equations containing memory terms that depend nonlinearly on past states model a variety of non-Markovian processes. In this study, we present a Markovian embedding procedure for such equations {with distributed delay} by utilising an exact spectral representation of the nonlinear memory function. This allows us to transform the nonlocal system to an equivalent local-in-time system in an abstract extended space. We demonstrate our embedding procedure for two  one-dimensional physical models: $(i)$ the walking droplet and $(ii)$ the single-phase Stefan problem. In addition to providing an alternative representation of the underlying physical system, the local representation finds applications in designing efficient time-integrators with time-independent computational costs for memory-dependent systems which typically suffer from growing-in-time costs.
\end{abstract}

\end{frontmatter}

\section{Introduction}{\label{sec1}}
Memory effects often emerge in reduced models of physical processes having many dynamically interacting degrees of freedom. They arise when one isolates a state variable of primary interest by integrating out the effects of other state variables. Mathematical models describing evolution of such primary state variables typically take the form of nonlocal differential equations with distributed delay. Examples include particle motion in unsteady hydrodynamic environments~\citep{lovalentiBrady93}, dynamics of self-propelled walking droplets~\citep{Oza2013} and chemically active particles~\citep{schnitzer23}, boundary evolution in diffusion processes in time-dependent domains~\citep{stefan1891,fokaspelloni2012} and under nonlinear boundary forcing~\citep{mannWolf51, nonlinearRadiation72, olmstead76}, dynamics of optical and mechanical resonators~\citep{Peters2022} and population dynamics~\citep{Cushing1977-rt}.

We are concerned with nonlocal models where the evolution equation for a state variable $y(t)$ has the following canonical form:
\begin{equation}{\label{eq:canonicalForm}}
    y^{(n)}(t) = L(t,y(t)) + \int_0^t N(y(t),y(s),t, s) \:ds~, 
\end{equation}
where the superscript $n$ indicates the order of time derivative. We distinguish between the local-in-time driving term, $L(\cdot)$, which depends only on the instantaneous state, and the nonlocal driving term, represented by the integral over history of states involving the integrable memory function $N(\cdot)$. The nonlocal term may be further classified based on whether the memory function depends linearly or nonlinearly on the history of states. The generalized Langevin equation \citep{Zwanzig1973} and the Maxey-Riley equation \citep{maxey1983equation,gatignol1983} are examples of the former. In this article, we consider the general case of a nonlinear memory function $N(\cdot)$.

The nonlocality poses challenges in solving equations such as \cref{eq:canonicalForm}.  In the absence of the nonlocal term, the dynamics is \textit{memory-independent}, and the resultant \cref{eq:canonicalForm} can be readily transformed into a system of ordinary differential equations (ODEs), yielding a Markovian description for which powerful tools from dynamical systems theory are available. However, with the nonlocal term intact, such a transformation is non-trivial. A common strategy is to realise Markovian embeddings in an abstract extended space by introducing a set of auxiliary variables that encode the memory effects, each governed by its own local evolution equation \citep{giorgi2010, HenrikNevermann_2023}. We show that a Markovian prescription can be naturally realised for nonlocal equations of the form in \cref{eq:canonicalForm} by an embedding procedure that relies on an exact spectral representation of the nonlinear memory function, extending ideas previously discussed for linear memory functions in \citep{Jaganathan2025}.
%\rv{The exploration of the more general case of Markovian embedding of equations with nonlinear memory kernels has been limited~\citep{HenrikNevermann_2023}.} 
%(\rvt{ For Divya: Is this study referening to only linear memory kernels?}). 
%\rvt{(For Divya: add any other exemplery studies here that look at Markovian embedding with nonlinear memory kernels if you know of any?).}

The paper is organized as follows. In section~\ref{sec1} we prescribe a general approach for Markovian embedding based on the spectral representation. Then, in \cref{sec2}, we demonstrate this embedding procedure to two different physical models, namely the 1D walking droplet in \cref{sec2p1} and the single-phase 1D Stefan problem in \cref{sec2p2} and validate our approach against known results in these systems. We conclude in section~\ref{sec4}.

\section{Spectral structure and Markovian embedding}\label{sec1}
Our approach involves expressing the nonlocal integral term in \cref{eq:canonicalForm} as a local-in-time term. Critical to the approach, we assume that the nonlinear memory function has a spectral representation of the following form:
\begin{equation}{\label{eq:genspectralform}}
    N(y(t),y(s),t, s) = \int_{\Gamma} e^{\phi(k;  y(t),y(s), t, s)} \psi(k,y(s)) \: dk,~~ 0< s < t~,
\end{equation}
where $\Gamma$ is a smooth contour in the complex plane and $\phi, \psi$ are  complex-analytic functions of the variable $k$. {Further, we assume $\phi$ to be additively separable into a local component, depending on the present state $(y(t),t)$, and a nonlocal component, depending on $(y(s), s)$, such that:
\[ \phi(k; y(t),y(s),t, s) = \phi_1(k;y(t),t) + \phi_2(k;y(s),s)~,\]
where each $\phi_i$, $i=\{1,2\}$, is complex-analytic as well.}

The proposed spectral representation allows embedding of \cref{eq:canonicalForm} into an extended space. This is done by substituting the spectral representation in the nonlocal term, followed by an interchange of order of integrals, to give:
\[\int_0^t \int_{\Gamma} e^{\phi(k; y(t), y(s), t, s)} \psi(k,y(s)) \: dk \:ds =: \int_{\Gamma} H(k,t) \:dk~, \]
\[\implies H(k,t) = \int_0^t e^{\phi(k;  y(t),y(s), t, s)} \psi(k,y(s))  \: ds~,\]
where $H(k,t)$ is the newly-introduced complex-valued auxiliary variable. Since it encodes the memory effect, we refer to it as the \textit{history function}. Owing to the particular spectral form in \cref{eq:genspectralform} and the additive-separability of $\phi$, the history function has a Markovian evolution given by an ODE, parameterised by the spectral variable $k$. Therefore, we have the following local-in-time reformulation of \cref{eq:canonicalForm} in an infinite-dimensional space:
\begin{subequations}{\label{eq:canonicalReformulation}}
    \begin{align}
        y^{(n)}(t) &= L(t,y(t)) + \int_{\Gamma} H(k,t) \:dk~, {\label{eq:canonicalReformulation_a}}\\
        \dot{H}(k,t) &=  \dot{\phi}(k;y(t),t) H(k,t) + e^{\phi(k;y(t),t)}\psi(k,y(t))~, {\label{eq:canonicalReformulation_b}} 
    \end{align}
\end{subequations}
where the overdot denotes {total} time derivative with respect to $t$. We emphasise that the reformulation is exact, and \cref{eq:canonicalReformulation} is entirely an equivalent representation of \cref{eq:canonicalForm}.

We note that this embedding, and hence the resulting history function, is not unique. The construction above is based on the proposed spectral representation. Other equivalent Markovian embeddings may however be realized.

% SECTION 2
\section{Prototypical physical models with nonlinear memory effects}{\label{sec2}}
We demonstrate our embedding procedure for two physical models with nonlinear memory effects, namely the one-dimensional walking droplet and the single-phase one-dimensional Stefan problem. Both models illustrate that a Markovian embedding in an infinite-dimensional space can be constructed subject to a natural spectral representation of the nonlinear memory kernel \eqref{eq:genspectralform}. These models differ in the complexity of the auxiliary history variable introduced upon embedding, highlighting the versatility of the approach.

% WALKER PROBLEM 
\subsection{Walking droplets}{\label{sec2p1}}
A hydrodynamic active system described by non-Markovian dynamics is that of self-propelled walking \citep{Couder2005WalkingDroplets} and superwalking~\citep{superwalker} droplets. By vertically vibrating an oil bath, a drop of the same oil can be made to bounce and walk on the liquid surface. Each bounce of the droplet locally excites a damped standing wave. The droplet interacts obliquely with these self-excited waves on subsequent bounces to propel itself horizontally, giving rise to a self-propelled, classical wave-particle entity (WPE). At large vibration amplitudes, the droplet-generated waves decay slowly in time. Hence, the motion of the droplet is affected by the history of waves along its trajectory. This gives rise to \textit{path memory} in the system and makes the dynamics non-Markovian. 
 
Oza \textit{et al.}~\citep{Oza2013} developed a theoretical stroboscopic model to describe the horizontal walking motion of such a WPE. The model averages over the fast vertical periodic bouncing of the droplet and provides a trajectory equation for the slow walking dynamics in the horizontal plane. We consider a reduction of this model to one horizontal dimension, $x\in\mathbb{R}$ (see Fig.~\ref{fig:typicalwalkerstates}a). Consider a droplet with position and velocity  given by $(x_d(t), \dot{x}_d(t)) \in \mathbb{R}^2$, which continuously generates standing waves with prescribed spatial structure $W(x)$ that are centered at the droplet position and decay in time. The dynamics of the $1$D WPE follows the non-dimensional integro-differential equation:
 \begin{equation}{\label{eq:non_markovianSWD}}
    \ddot{x}_d(t)= - \dot{x}_d(t) - C_1\int_0^t W'(x_d(t)-x_d(s))K(t-s;C_2)\: ds, 
\end{equation}
where $C_1$ and $C_2$ are non-negative constants representing dimensionless wave-amplitude and inverse memory parameter, respectively\footnote{Note that $C_1$ and $C_2$ are related to the dimensionless parameters $\kappa$ and $\beta$ in Oza \textit{et al.}~\citep{Oza2013} by $C_1=\beta \kappa^2$ and $C_2=\kappa$. For more details see \citep{Valani2024infmem}.}. We refer the reader to Ref.~\citep{Oza2013} for details and explicit expressions for these parameters.  \Cref{eq:non_markovianSWD} is a horizontal force balance of the WPE with the right-hand side containing an effective drag term proportional to velocity $-\dot{x}_d(t)$ and the nonlocal memory term capturing the cumulative force on the particle from the superposition of the self-generated waves along its path. The memory kernel comprises the functions $W'(\cdot)$ and $K(\cdot)$;  the former represents the wave-gradient where the prime denotes derivative with respect to its argument, and $K(\cdot)$ imposes the temporal decay. In the stroboscopic model of a walking droplet, $-W'(x) = J_1(x)$, where $J_1$ is the Bessel-$J$  function of order one and $K(t)=e^{-C_2t}$.

In the high-memory regime ($C_2\ll 1$), WPEs exhibit hydrodynamic quantum analogs~\citep{Bush_2021}. However, the regime may become experimentally difficult to access~\citep{PhysRevLett.122.104303} due to the increased susceptibility of the system to the Faraday instability~\citep{Faraday1831a}. Numerical simulations provide an alternative with controllability to probe this regime, especially in relation to investigating hydrodynamic quantum analogs, but they also entail dealing with the non-Markovian structure of \cref{eq:non_markovianSWD} and time-dependent computational costs therein. 
\subsubsection{Markovian embedding for the stroboscopic model}
We convert Eq.\eqref{eq:non_markovianSWD} to a Markovian description in the following way. We recall the following integral representation of the Bessel-$J_1$ function for some $z\in \mathbb{R}$:
\begin{equation}
    J_1(z) = -\frac{i}{\pi}\int_{-1}^1 e^{ikz} \frac{k}{\sqrt{1-k^2}} \:dk ~.
\end{equation}
Substituting the above in the memory term of \cref{eq:non_markovianSWD}, followed by a switch in the order of integrals, we construct the equivalent local-in-time representation for the memory integral,
\begin{equation*}
-\frac{i}{\pi} \int_0^t \int_{-1}^1 e^{ik(x_d(t)-x_d(s))-C_2(t-s)} \frac{k}{\sqrt{1-k^2}} \:dk \:ds  =: \int_{-1}^1 H(k,t) w(k) \:dk  
\end{equation*}
where the weight function $w(k) = 1/\sqrt{1-k^2}$ and $H(k,t)$ is a complex-valued function of time $t$ and a real number $k$ with a finite support in $[-1,1]$. The induced definition of $H(k,t)$ is
\begin{equation}{\label{Hdef}}
H(k,t) := -\frac{i k}{\pi} \int_0^t e^{ik(x_d(t)-x_d(s))-C_2(t-s)} \:ds~.
\end{equation}
The form of $H(\cdot)$ in \eqref{Hdef} suggests that it has a Markovian evolution according to an ODE parameterised by the spectral variable $k$. Consequently, combined with the definition of the droplet's velocity $\dot{x}_d = v_d$, we derive the following Markovian prescription for the WPE dynamics in the extended state space for $t>0$:
\begin{subequations}{\label{eq:markovianSWD}}
    \begin{align}
        \dot{v}_d(t) &= -v_d(t) + C_1 \int_{-1}^{1} H(k,t)w(k)  \:dk, \label{eq:3a}\\
        \dot{H}(k,t) &= -C_2 H(k,t) + ik v_d(t) H(k,t) -\frac{ik}{\pi} \label{eq:3b}~,
    \end{align}
\end{subequations}
subject to initial conditions $(x_{d0}, v_{d0})$ and $H(k,0)=0$. We note in \cref{Hdef} that $H(\cdot)$ preserves certain symmetries with respect to the spectral variable $k$ at all times: $\text{Re}(H)$ has an even-symmetry whereas $\text{Im}(H)$ is odd-symmetric. Therefore, whereas the real and imaginary parts of the history function drive each other's dynamics, only the real part contributes to the memory integral in \cref{eq:3a}.

\begin{figure}[h!]
    \centering
    \subfloat[]{{\includegraphics[width=0.237\textwidth]{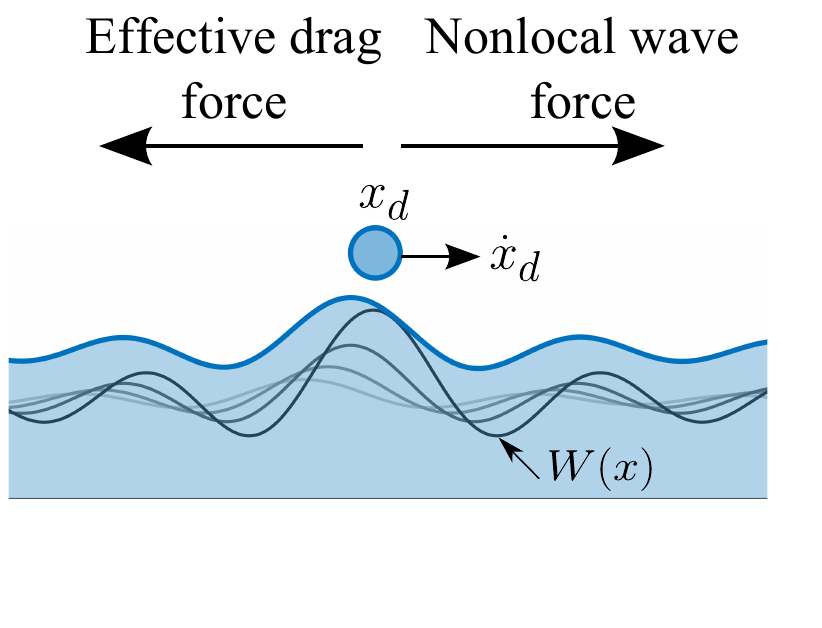}}}
    \:
    \subfloat[]{{\includegraphics[width=0.235\textwidth]{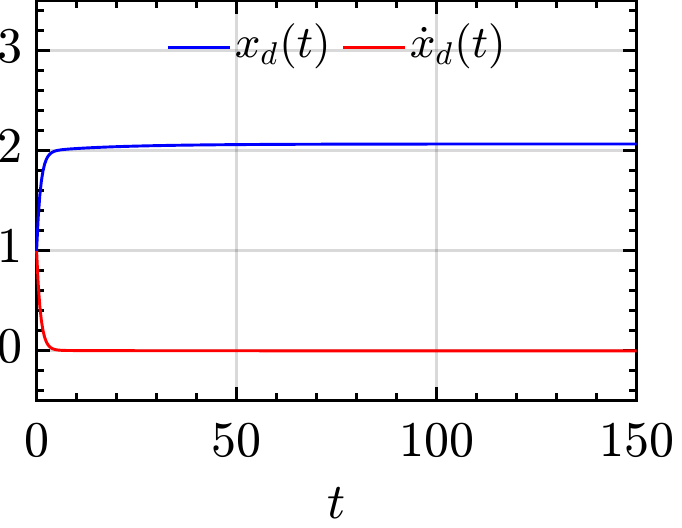}}}%
    \!
    \subfloat[]{{\includegraphics[width=0.245\textwidth]{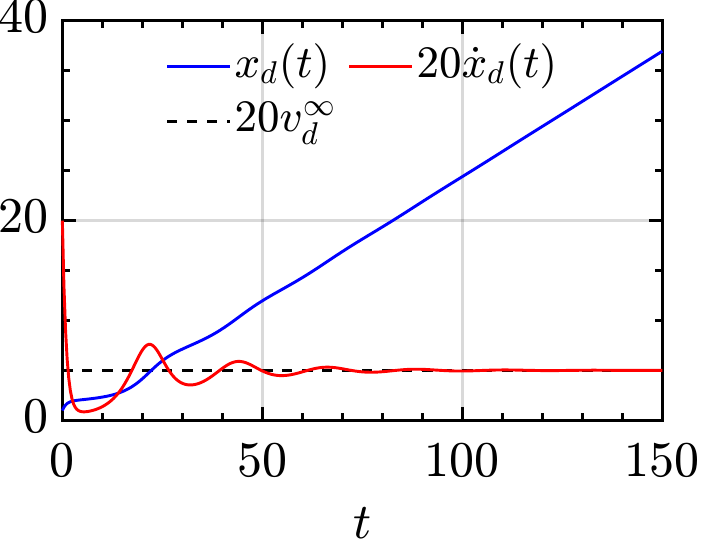}}}%
    \!
    \subfloat[]{{\includegraphics[width=0.245\textwidth]{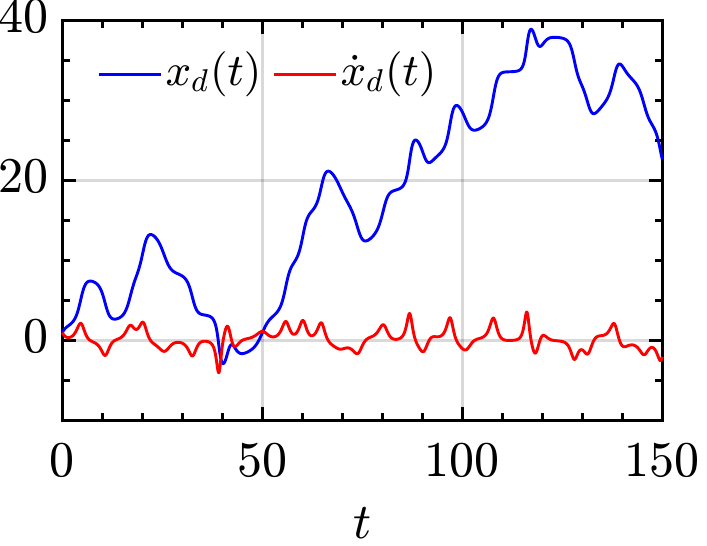}}}%
    \captionsetup{width=\linewidth, font=footnotesize}
    \caption{(a) Schematic of a 1D walking droplet. Typical known droplet states in the stroboscopic model obtained by solving the Markovian system (\cref{eq:markovianSWD}) for $(x_{d0}, \dot{x}_{d0})$ = $(1,1)$: (b) \textit{Non-walker} $(C_1=0.01, C_2=0.1)$, (c) \textit{Steady walker }$(C_1=C_2=0.1)$, (d) \textit{Chaotic walker }$(C_1=1.5,C_2=0.01)$. Velocity in (c) is scaled by a factor of $20$ for visibility.}%
    \label{fig:typicalwalkerstates}
\end{figure}

\begin{figure}[h!]
    \captionsetup{width=\linewidth, font=footnotesize}
    \subfloat{{\includegraphics[width=0.25\textwidth]{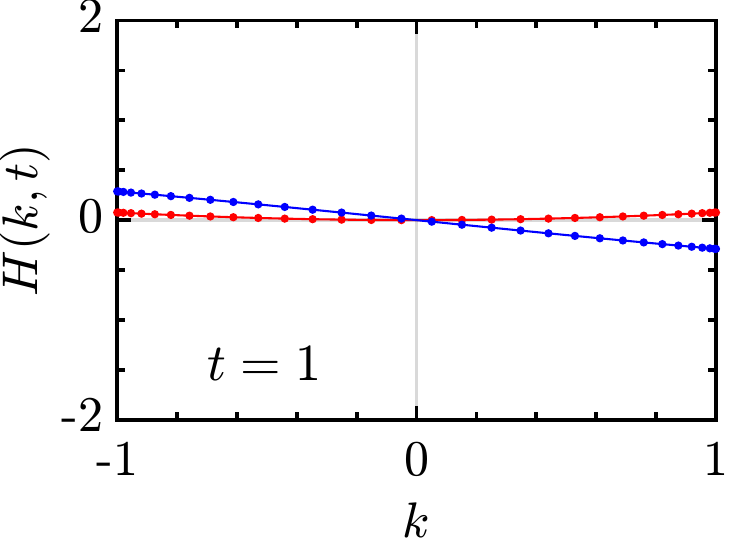} }}%
    \!
    \subfloat{{\includegraphics[width=0.23\textwidth]{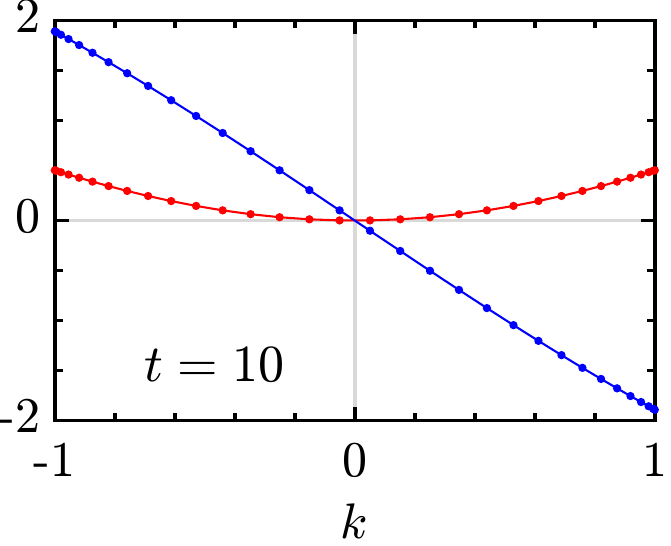} }}%
    \!
    \subfloat{{\includegraphics[width=0.23\textwidth]{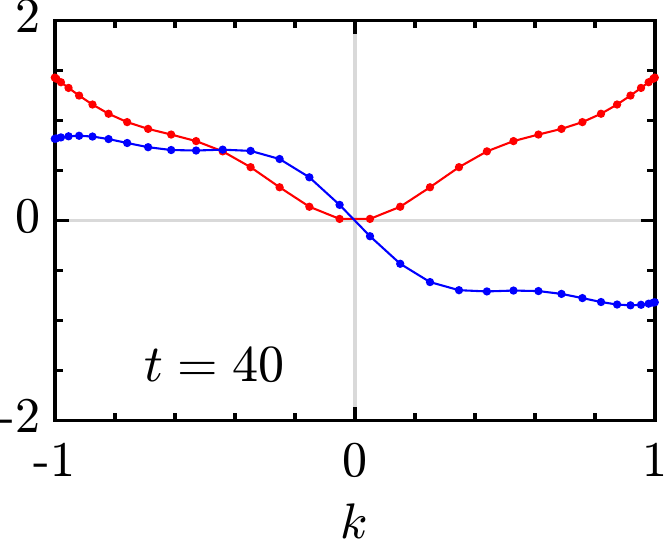} }}%
    \!
     \subfloat{{\includegraphics[width=0.23\textwidth]{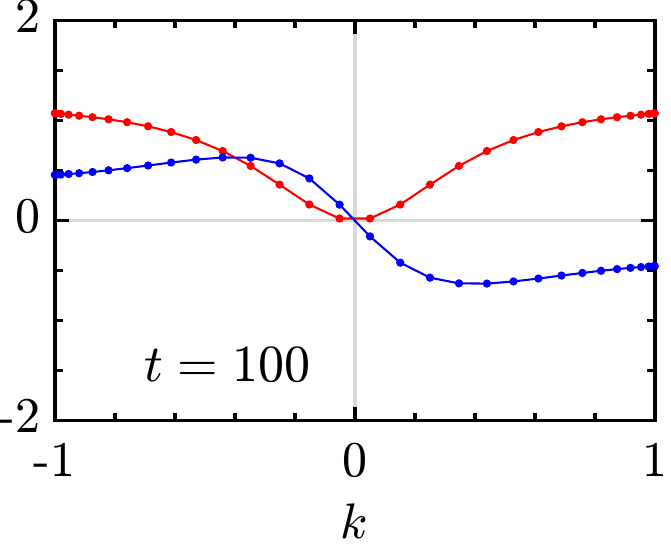} }}%
    \caption{Evolution of real (red)/imaginary (blue) parts of the history function $H(k,t)$ in the spectral space at representative times for a steady walker $(C_1=C_2=0.1)$, with $H(k,0)=0$. The finite support of $H(k,t)$ in $[-1,1]$ and its smoothness demands a nominal, fixed (in time) requirement of as few as $M=30$ Chebyshev quadrature nodes to accurately compute the history integral.}
    \label{fig:walkerHevolution}% 
\end{figure}

The resultant set of local differential equations~\eqref{eq:markovianSWD} can be readily solved using any standard time-integrator; we use the second-order Runge-Kutta scheme. An additional task involves computing the history integral over $k$. The integrand, with its finite support in $[-1,1]$ and the form of weight function $w$, naturally suggests expansion of $H(k,t)$ in the bases of Chebyshev polynomials of the first kind. Therefore, we use the spectrally-accurate Clenshaw-Curtis quadrature method to approximate the integral:
\[\int_{-1}^1 H(k,t) w(k) \:dk \approx \sum_{n=0}^M\omega_n H(k_n,t), \quad M \in \mathbb{N}~,\]
where $k_n = \cos(n\pi/M)$ are the Chebyshev nodes and $\omega_n$ are the associated weights. In order to validate the Markovian prescription derived for the WPE dynamics, we numerically solve \cref{eq:markovianSWD} for a few representative parameter sets $(C_1,C_2)$. Fig.~\ref{fig:typicalwalkerstates} shows that the embedded system of equations~\eqref{eq:markovianSWD} successfully reproduces the previously known non-walking and walking regimes in the parameter space~\citep{Durey2020,Valaniunsteady2021}. Further, for a droplet in the steady walking regime, an analytical expression for its steady speed~\citep{Oza2013} is :
\begin{equation}{\label{eq:steadyspeed}}
    v_d^{\infty} = \frac{1}{\sqrt{2}}\sqrt{2C_1-C_2^2-\sqrt{C_2^4+4C_1C_2^2}}~.
\end{equation}
The numerical solution for the steady walker attains the above analytical steady walking speed (dashed line) in Fig.~\ref{fig:typicalwalkerstates}c. Additionally, in Fig.~\ref{fig:walkerHevolution}, we show the evolution of the newly-introduced auxiliary history function, $H(\cdot)$, in the finite $k-$domain over time. The spectral behaviour of $H(\cdot)$ suggests that a low-dimensional quadrature approximation with fewer Chebyshev nodes, $M$, is adequate to accurately compute its integral over the finite spectral domain.

There have been previous works~\citep{phdthesismolacek,Durey2020,Durey2020lorenz,Valanilorenz2022, Valani2024infmem} that rewrite the integro-differential equation for the WPE into a system of ODEs. However, these transformations are specific to the choices of the wave form $W(x)$. For instance, for a similar choice of $W' = -J_1$, \citep{Durey2020} use the recurrence relations for the Bessel-$J$ function to derive a system of ODEs. The Markovian embedding formalism is applicable for a broader class of wave forms that have a suitable spectral representation. This is particularly useful in generalised pilot-wave framework, where new hydrodynamic quantum analogues are being explored by investigating various wave forms~\citep{Bush_2021}.

% STEFAN PROBLEM
\subsection{Single-phase Stefan problem}{\label{sec2p2}}
We now consider the well-known class of free boundary problems called the \textit{Stefan problem}, which primarily describes phase-change processes such as the melting of a solid \citep{stefan1891, guentherlee2012}. In its simplest non-dimensional formulation, it comprises a one-dimensional domain in $\mathbb{R}^+$, contiguously supporting a molten phase and a solid phase, separated at their interface, which is free to move as the solid melts (see Fig.~\ref{fig:stefansol}a). The solid phase is modelled as an infinite heat sink maintained at the melting temperature at all times. Therefore, the simplified problem involves finding the solution pair $(\theta(x,t), l(t))$, where $\theta(x,t)$  describes the instantaneous temperature distribution in the molten phase and $l(t)$ is the location of the melting front. The function $\theta$ satisfies the diffusion equation $\partial_t \theta - \partial_x^2 \theta = 0$ in $x\in [0,l(t)]$, subject to an arbitrary initial condition $\theta_0(x)$ in the initial domain $x \in [0,l_0]$ at $t=t_0$, and a temperature or heat flux condition at the fixed boundary $x=0$. The moving front, which is at the melting temperature, is governed by the \textit{Stefan condition} $\dot{l} = -\partial_x \theta(l(t),t)$. Here, the notation $\partial_{(\cdot)}$ denotes a partial derivative with respect to the subscript variable.

We consider the case where temperature is prescribed at the fixed boundary, $\theta(0,t)=f(t)$, for exposition. With primary interest in the dynamics of the interface's location, the bulk heat diffusion process in the molten phase may be effectively ``integrated out'' to derive a non-Markovian equation of motion for the moving front. The resulting velocity equation for the moving front, $v(t) = \dot{l}(t)$, is compactly written in the following nonlinear Volterra integral form for $t>t_0$ \citep{fokaspelloni2012,guentherlee2012}:
\begin{equation}\label{eq:movingfront}
v(t) = g(t-t_0,l(t);\theta_0', l_0) + \int_{t_0}^t N(l(t),l(s), t-s; v(s), \dot{f}(s)) \:ds
\end{equation}
where $\theta_0', \dot{f}$ denote the  spatial derivative and temporal derivative of $\theta_0, f$ respectively, and the function $g$ is:
\begin{equation}
     g(t-t_0,l(t);\theta_0', l_0) = -\frac{1}{\sqrt{\pi (t-t_0)}}\int_0^{l_0} \Big( e^{-(l(t)+x)^2/4(t-t_0)}+e^{-(l(t)-x)^2/4(t-t_0)}\Big) \theta_0'(x) \:dx~.
\end{equation}
The nonlinear kernel may be decomposed into contributions from forcing at the fixed boundary and the unknown velocity of the solid-liquid interface as follows:
\begin{equation}{\label{def:nonlinearKernel}}
    N(l(t), l(s), t-s;v(s),\dot{f}(s)) =  v(s) N_1(l(t),l(s),t-s) + \dot{f}(s)N_2(l(t),0,t-s)~,
\end{equation} 
with the definitions
\begin{equation*}
        N_1(x,y,z) =\frac{1}{2\sqrt{\pi}}\Bigg(\frac{(x+y)e^{-(x+y)^2/4z} - (x-y)e^{-(x-y)^2/4z}}{z^{3/2}} \Bigg), \quad N_2(x,y,z) = \frac{2}{\sqrt{\pi}}\frac{e^{-x^2/4z}}{z^{1/2}}~.
\end{equation*}
Note that the function $g$ is a local-in-time term. The second term on the right-hand side in \cref{eq:movingfront} is the memory term, which introduces non-locality and the nonlinear dependence on the moving front $l(t)$.

\subsubsection{Markovian embedding for Stefan problem}
As before, we construct an embedding such that the present non-Markovian representation for $v(t)$ may be turned Markovian. We identify the following exact spectral representation of the nonlinear kernel for a real $k$:
\begin{subequations}
\begin{align}
N_1(x,y,z) &= \int_{-\infty}^{\infty} n_1(k;x,y,z) \:dk = \frac{i}{\pi}\int_{-\infty}^{\infty}  \Big(e^{-k^2z + ik(x-y)} -e^{-k^2z + ik(x+y)} \Big) k \:dk, \label{eq:7a}\\
N_2(x,y,z) &= \int_{-\infty}^{\infty} n_2(k;x,y,z) \: dk= \frac{1}{\pi}\int_{-\infty}^{\infty} \Big(e^{-k^2z+ik(x-y)}+e^{-k^2z+ik(x+y)}\Big) \:dk~.
\end{align}
\end{subequations}
Substituting the above spectral representations in the memory term, followed by a switch in the order of integrals, we derive the local representation with the introduction of the auxiliary history function $H(k,t)$:
\begin{equation*}
    \int_{t_0}^t \Bigg(v(s)\int_{-\infty}^{\infty} n_1(k;l(t),l(s),t-s)\:dk + \dot{f}(s)\int_{-\infty}^{\infty}n_2(k;l(t),0,t-s)\:dk \Bigg)\:ds =: \int_{-\infty}^{\infty} H(k,t) \:dk~.
\end{equation*}
The corresponding induced definition of the complex-valued history function is:
\begin{equation}\label{history4stefan}
    H(k,t) := \int_{t_0}^t \Big(v(s) n_1(k;l(t),l(s),t-s) + \dot{f}(s) n_2(k;l(t),0,t-s)\Big) \:ds  ~.
\end{equation}
Differentiating the above with respect to time, one may derive an ODE for the history function and realise the following equivalent Markovian prescription for the moving front for $t>t_0$:
\begin{subequations}{\label{eq:markovianStefan}}
    \begin{align}
         \dot{l}(t) &= g(t-t_0,l(t);\theta_0', l_0) + \int_{-\infty}^{\infty} H(k,t) \: dk, \label{eq:markovianStefan_1}\\
        \dot{H}(k,t)  &= -k^2 H(k,t) +ikv(t) H(k,t) + \frac{ik}{\pi} (1-e^{2ikl(t)}) v(t) + \frac{2}{\pi} e^{ikl(t)} \dot{f}(t) \label{eq:markovianStefan_2} ~,
    \end{align}
\end{subequations}
subject to $l(t_0)=l_0, v(t_0)=g(0,l_0;\theta_0', l_0), \; H(k,t_0) = 0$. The history function in this case too preserves similar symmetries, suggesting that only its even-symmetric real part contributes to the history integral.

We show equivalence of the derived embedded Markovian system to the original non-Markovian system (\cref{eq:movingfront}) by numerically solving \cref{eq:markovianStefan}. We use the second-order Runge-Kutta exponential time-differencing method \citep{CM2002} to solve for $H(k,t)$ due to stiffness introduced by the $-k^2 H$ term, and a standard second-order time-integrator to solve for $l(t)$. The latter requires evaluating the history integral whose quadrature approximation, however, demands different treatment from the walker problem on two accounts:
 \begin{itemize}[leftmargin=*]
        \item[1.] $H(k,t)$ has infinite support in the $k-$space. While this warrants truncation of the $k-$space, its decay behaviour at large $k$ constraints the extent of truncation.
        \item[2.] $H(k,t)$ is highly oscillatory; the frequency of oscillations increases with both $k$ and $t$, which is ascribed to terms such as $e^{ikl(t)}$ in \cref{eq:markovianStefan_2}. The dependence on $t$ through $l(t)$ exacerbates the oscillations in domain growth problems such as the one under discussion. Consequently, for an accurate quadrature approximation, an increasingly dense set of collocation points in the truncated domain is required.
\end{itemize}
The above points are cautionary observations. While one could potentially address these concerns through computationally efficient methods, such an undertaking exceeds the scope of our present work. Therefore, for our purpose of demonstrating equivalence, we adopt a heuristic approach to compute the history integral. This involves truncating the $k-$space, mapping it to the interval $[-1,1]$, and employing Clenshaw-Curtis quadrature to compute the history integral.

\begin{figure}[h!]
    \centering    
    \subfloat[]{{\includegraphics[width=0.37\textwidth]{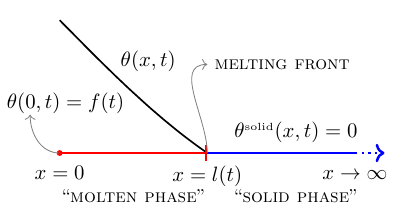} }}%
    \quad
    \subfloat[]{{\includegraphics[width=0.24\textwidth]{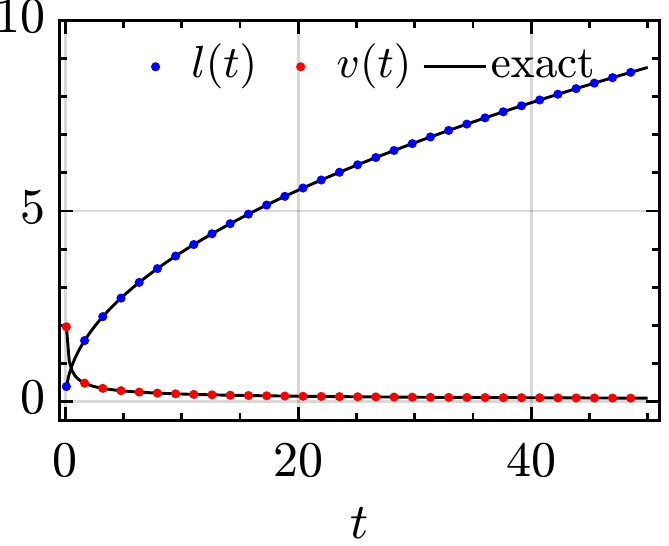} }}%
    \quad
    \subfloat[]{{\includegraphics[width=0.28\textwidth]{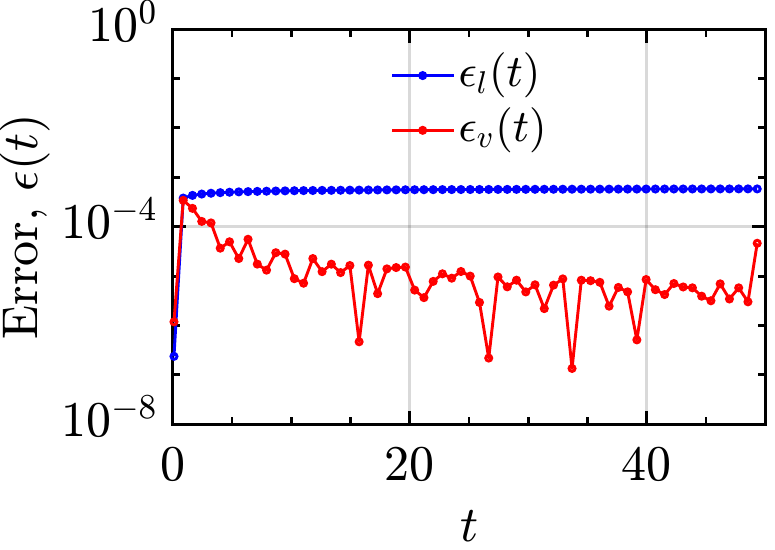} }}%
    \captionsetup{width=\linewidth, font=footnotesize}
    \caption{(a) Schematic of 1D single-phase Stefan problem. (b) Numerical response of the melting front to a constant temperature forcing $\theta(0,t)=f(t)=1$ at the fixed end $x=0$ overlaid on the exact solution along with (c) the instantaneous pointwise error.}
    \label{fig:stefansol}
\end{figure}

\begin{figure}[h!]
    \centering
    \captionsetup{width=\linewidth, font=footnotesize}
    \subfloat{{\includegraphics[width=0.3\textwidth]{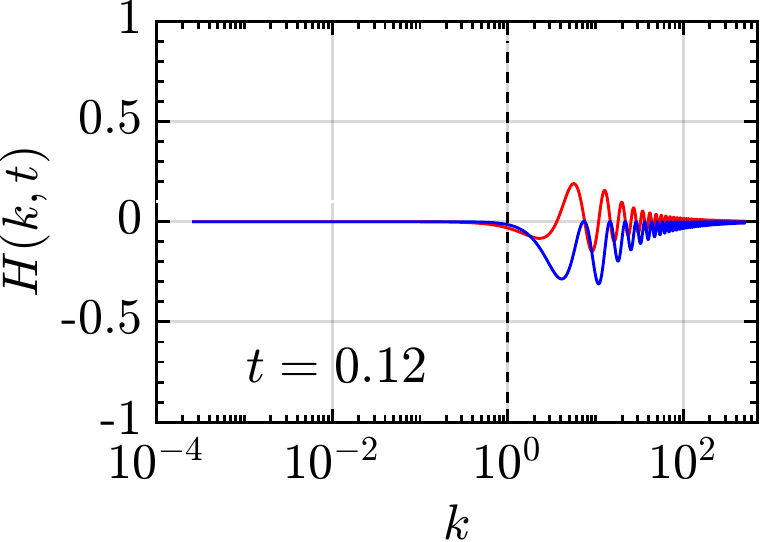} }}%
    \qquad
    \subfloat{{\includegraphics[width=0.28\textwidth]{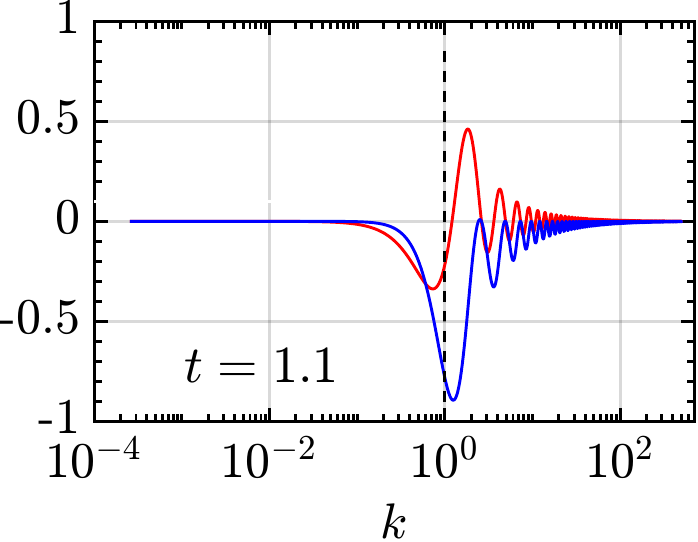} }}%
    \qquad
     \subfloat{{\includegraphics[width=0.28\textwidth]{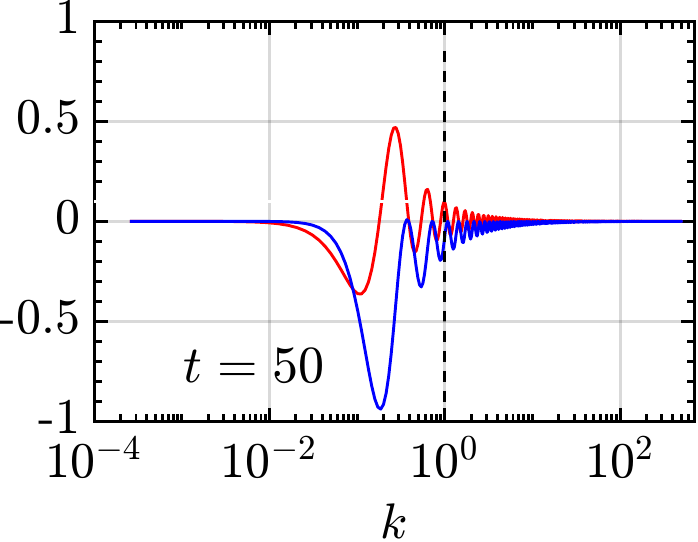} }}%
    \caption{Evolution of real (red)/imaginary (blue) parts of the history function $H(k,t)$ in the truncated $k-$domain $k \in [0,500]$ at different time instances for the single-phase Stefan problem subject to constant temperature forcing $(f(t)=1)$ at the fixed end $x=0$.  $M=2000$ Chebyshev nodes were used to accurately compute the integral of the highly oscillatory history function.} 
    \label{fig:stefanHevolution}%
\end{figure}

We consider the example corresponding to melting due to constant temperature at the fixed end, $f(t) = 1$, with the following analytical solution pair \citep{mitchell2009}:
\begin{equation}{\label{eq:exactinetrfaceloc}}
    \theta(x,t) =  1 - \frac{\text{erf}(x/2\sqrt{t})}{\text{erf}(\alpha)} \quad \text{for } x \in [0,l(t)], \quad l(t) = 2\alpha \sqrt{t}, \quad t>0~,
\end{equation}
where the constant $\alpha$ satisfies the transcendental equation: $\sqrt{\pi} \alpha \exp(\alpha^2) \text{erf}(\alpha) = 1$ and $\text{erf}(\cdot)$ is the error function. To avoid the degeneracy at $t=0$ due to zero-length domain in this example, we let the process evolve for time $t_0>0$ to a non-zero domain length $l_0$. Prescribing $(\theta(x,t_0), l_0)$ as the initial state, we numerically evolve \cref{eq:markovianStefan} from $t_0$. Fig.~\ref{fig:stefansol} shows agreement between the numerical and analytical solutions (\cref{eq:exactinetrfaceloc}) for location of the moving front, supplemented with the pointwise error. In Fig.~\ref{fig:stefanHevolution}, we plot the pertinent history function in the truncated spectral space at different time instances. The highly oscillatory behaviour of $H(\cdot)$ in the truncated $k-$domain is evident. Consequently, we are forced to include a larger number of Chebyshev nodes to numerically integrate the history function.

% Section 3
\section{Conclusions}\label{sec4}
We have described a Markovian embedding procedure for a class of  evolutionary equations with memory effects, which critically relied on the exact spectral representation of the memory kernel. We have explicitly shown the embedding procedure for two physical models in Section~\ref{sec2}, namely the one-dimensional walking droplet and the single-phase Stefan problem. Notably, in both cases, the memory kernel is a nonlinear function of the underlying state variable and has an exact spectral representation. We have validated the resultant Markovian prescriptions for each model by reproducing known results of the respective underlying systems.

Physical processes inherently follow Markovian dynamics when described adequately by all the driving state variables. The non-Markovian descriptions, such as in Eqs.~\eqref{eq:non_markovianSWD}, \eqref{eq:movingfront}, are often the result of isolating the dynamics of a state variable by ``integrating out" the effects of ``environment" comprising other state variables. While identifying these integrated physical variables may not always be feasible, our Markovian embedding procedure provides an alternative mathematical reconstruction of the Markovian dynamics. 

From a computational standpoint, a Markovian representation ensures that the numerical evolution of the corresponding time-discretized system incurs a time-independent cost. This is in contrast with the standard approaches for memory-dependent systems, where the computational cost grows with time. It is important to recognise that our Markovian embedding procedure comes at the cost of solving an additional local-in-time equation for the auxiliary history function $H(\cdot,t)$, which is an infinite-dimensional object. However, the associated computational cost is independent of time, and relies on the finite-dimensional reduction of the history function in the spectral space. An accurate finite-dimensional approximation depends on the behaviour of $H(\cdot,t)$ in the spectral space. In this regard, the two model problems that we have discussed demonstrate the extreme scenarios: the Stefan problem required a higher-dimensional approximation of the history function, comprising thousands of spectral variables, while the walker problem allowed a lower-dimensional approximation, with only a few tens of spectral variables for an accurate representation of the history function.

% Acknowledgements
\section*{Acknowledgements}
DJ acknowledges support of the Department of Atomic Energy, Government of India, under project no. RTI4001. RV was supported by Australian Research Council (ARC) Discovery Project DP200100834 during the course of the work. We thank Vishal Vasan for introducing the Stefan problem to us and for the various discussions pertaining to this work. We also thank Rama Govindarajan for her feedback on the draft. 

\bibliographystyle{model1-num-names}
\bibliography{references}

\begin{thebibliography}{30}
\expandafter\ifx\csname natexlab\endcsname\relax\def\natexlab#1{#1}\fi
\providecommand{\bibinfo}[2]{#2}
\ifx\xfnm\relax \def\xfnm[#1]{\unskip,\space#1}\fi
%Type = Article
\bibitem[{Lovalenti and Brady(1993)}]{lovalentiBrady93}
\bibinfo{author}{P.~M. Lovalenti}, \bibinfo{author}{J.~F. Brady},
\newblock \bibinfo{title}{The hydrodynamic force on a rigid particle undergoing
  arbitrary time-dependent motion at small {R}eynolds number},
\newblock \bibinfo{journal}{J. Fluid Mech.} \bibinfo{volume}{256}
  (\bibinfo{year}{1993}) \bibinfo{pages}{561–605}.
%Type = Article
\bibitem[{Oza et~al.(2013)Oza, Rosales, and Bush}]{Oza2013}
\bibinfo{author}{A.~U. Oza}, \bibinfo{author}{R.~R. Rosales},
  \bibinfo{author}{J.~W.~M. Bush},
\newblock \bibinfo{title}{A trajectory equation for walking droplets:
  hydrodynamic pilot-wave theory},
\newblock \bibinfo{journal}{J. Fluid Mech.} \bibinfo{volume}{737}
  (\bibinfo{year}{2013}) \bibinfo{pages}{552--570}.
%Type = Article
\bibitem[{Peng and Schnitzer(2023)}]{schnitzer23}
\bibinfo{author}{G.~G. Peng}, \bibinfo{author}{O.~Schnitzer},
\newblock \bibinfo{title}{Weakly nonlinear dynamics of a chemically active
  particle near the threshold for spontaneous motion ii. history-dependent
  motion},
\newblock \bibinfo{journal}{Phys. Rev. Fluids} \bibinfo{volume}{8}
  (\bibinfo{year}{2023}) \bibinfo{pages}{033602}.
%Type = Article
\bibitem[{Stefan(1891)}]{stefan1891}
\bibinfo{author}{J.~Stefan},
\newblock \bibinfo{title}{Ueber die theorie der eisbildung, insbesondere über
  die eisbildung im polarmeere},
\newblock \bibinfo{journal}{Annalen der Physik} \bibinfo{volume}{278}
  (\bibinfo{year}{1891}) \bibinfo{pages}{269--286}.
%Type = Article
\bibitem[{Fokas and Pelloni(2012)}]{fokaspelloni2012}
\bibinfo{author}{A.~S. Fokas}, \bibinfo{author}{B.~Pelloni},
\newblock \bibinfo{title}{Generalized {D}irichlet-to-{N}eumann map in
  time-dependent domains},
\newblock \bibinfo{journal}{Stud. Appl. Math.} \bibinfo{volume}{129}
  (\bibinfo{year}{2012}) \bibinfo{pages}{51--90}.
%Type = Article
\bibitem[{Mann and Wolf(1951)}]{mannWolf51}
\bibinfo{author}{W.~R. Mann}, \bibinfo{author}{F.~Wolf},
\newblock \bibinfo{title}{Heat transfer between solids and gasses under
  nonlinear boundary conditions},
\newblock \bibinfo{journal}{Q. Appl. Math.} \bibinfo{volume}{9}
  (\bibinfo{year}{1951}) \bibinfo{pages}{163--184}.
%Type = Article
\bibitem[{Keller and Olmstead(1972)}]{nonlinearRadiation72}
\bibinfo{author}{J.~B. Keller}, \bibinfo{author}{W.~E. Olmstead},
\newblock \bibinfo{title}{Temperature of a nonlinearly radiating semi-infinite
  solid},
\newblock \bibinfo{journal}{Q. Appl. Math.} \bibinfo{volume}{29}
  (\bibinfo{year}{1972}) \bibinfo{pages}{559--566}.
%Type = Article
\bibitem[{Olmstead and Handelsman(1976)}]{olmstead76}
\bibinfo{author}{W.~E. Olmstead}, \bibinfo{author}{R.~A. Handelsman},
\newblock \bibinfo{title}{Asymptotic solution to a class of nonlinear
  {V}olterra integral equations. {II}},
\newblock \bibinfo{journal}{SIAM J. Appl. Math.} \bibinfo{volume}{30}
  (\bibinfo{year}{1976}) \bibinfo{pages}{180--189}.
%Type = Article
\bibitem[{Peters and Rodriguez(2022)}]{Peters2022}
\bibinfo{author}{K.~J.~H. Peters}, \bibinfo{author}{S.~R.~K. Rodriguez},
\newblock \bibinfo{title}{Limit cycles and chaos induced by a nonlinearity with
  memory},
\newblock \bibinfo{journal}{Eur. Phys. J. Spec. Top.} \bibinfo{volume}{231}
  (\bibinfo{year}{2022}) \bibinfo{pages}{247--254}.
%Type = Book
\bibitem[{Cushing(1977)}]{Cushing1977-rt}
\bibinfo{author}{J.~M. Cushing}, \bibinfo{title}{Integrodifferential equations
  and delay models in population dynamics}, Lecture Notes in Biomathematics,
  \bibinfo{publisher}{Springer}, \bibinfo{address}{Berlin, Germany},
  \bibinfo{year}{1977}.
%Type = Article
\bibitem[{Zwanzig(1973)}]{Zwanzig1973}
\bibinfo{author}{R.~Zwanzig},
\newblock \bibinfo{title}{Nonlinear generalized {L}angevin equations},
\newblock \bibinfo{journal}{J. Stat. Phys.} \bibinfo{volume}{9}
  (\bibinfo{year}{1973}) \bibinfo{pages}{215--220}.
%Type = Article
\bibitem[{Maxey and Riley(1983)}]{maxey1983equation}
\bibinfo{author}{M.~R. Maxey}, \bibinfo{author}{J.~J. Riley},
\newblock \bibinfo{title}{Equation of motion for a small rigid sphere in a
  nonuniform flow},
\newblock \bibinfo{journal}{Phys. Fluids} \bibinfo{volume}{26}
  (\bibinfo{year}{1983}) \bibinfo{pages}{883--889}.
%Type = Article
\bibitem[{Gatignol(1983)}]{gatignol1983}
\bibinfo{author}{R.~Gatignol},
\newblock \bibinfo{title}{The {F}ax{\'{e}}n formulae for a rigid particle in an
  unsteady non-uniform {S}tokes flow},
\newblock \bibinfo{journal}{J. Mec. Theor. Appl.} \bibinfo{volume}{2}
  (\bibinfo{year}{1983}) \bibinfo{pages}{241--282}.
%Type = Article
\bibitem[{Fabrizio et~al.(2010)Fabrizio, Giorgi, and Pata}]{giorgi2010}
\bibinfo{author}{M.~Fabrizio}, \bibinfo{author}{C.~Giorgi},
  \bibinfo{author}{V.~Pata},
\newblock \bibinfo{title}{A new approach to equations with memory},
\newblock \bibinfo{journal}{Arch. Ration. Mech. Anal.} \bibinfo{volume}{198}
  (\bibinfo{year}{2010}) \bibinfo{pages}{189--232}.
%Type = Article
\bibitem[{Nevermann and Gros(2023)}]{HenrikNevermann_2023}
\bibinfo{author}{D.~H. Nevermann}, \bibinfo{author}{C.~Gros},
\newblock \bibinfo{title}{Mapping dynamical systems with distributed time
  delays to sets of ordinary differential equations},
\newblock \bibinfo{journal}{J. Phys. A: Math. Theor.} \bibinfo{volume}{56}
  (\bibinfo{year}{2023}) \bibinfo{pages}{345702}.
%Type = Article
\bibitem[{Jaganathan et~al.(2025)Jaganathan, Govindarajan, and
  Vasan}]{Jaganathan2025}
\bibinfo{author}{D.~Jaganathan}, \bibinfo{author}{R.~Govindarajan},
  \bibinfo{author}{V.~Vasan},
\newblock \bibinfo{title}{Explicit integrators for nonlocal equations: The case
  of the {M}axey-{R}iley-{G}atignol equation},
\newblock \bibinfo{journal}{Quart. Appl. Math.} \bibinfo{volume}{83}
  (\bibinfo{year}{2025}) \bibinfo{pages}{135--158}.
%Type = Article
\bibitem[{Couder et~al.(2005)Couder, Proti{\`{e}}re, Fort, and
  Boudaoud}]{Couder2005WalkingDroplets}
\bibinfo{author}{Y.~Couder}, \bibinfo{author}{S.~Proti{\`{e}}re},
  \bibinfo{author}{E.~Fort}, \bibinfo{author}{A.~Boudaoud},
\newblock \bibinfo{title}{Dynamical phenomena: Walking and orbiting droplets},
\newblock \bibinfo{journal}{Nature} \bibinfo{volume}{437}
  (\bibinfo{year}{2005}) \bibinfo{pages}{208--208}.
%Type = Article
\bibitem[{Valani et~al.(2019)Valani, Slim, and Simula}]{superwalker}
\bibinfo{author}{R.~N. Valani}, \bibinfo{author}{A.~C. Slim},
  \bibinfo{author}{T.~Simula},
\newblock \bibinfo{title}{Superwalking droplets},
\newblock \bibinfo{journal}{Phys. Rev. Lett} \bibinfo{volume}{123}
  (\bibinfo{year}{2019}) \bibinfo{pages}{024503}.
%Type = Article
\bibitem[{Valani(2024)}]{Valani2024infmem}
\bibinfo{author}{R.~N. Valani},
\newblock \bibinfo{title}{{Infinite-memory classical wave-particle entities,
  attractor-driven active particles, and the diffusionless Lorenz equations}},
\newblock \bibinfo{journal}{Chaos} \bibinfo{volume}{34} (\bibinfo{year}{2024})
  \bibinfo{pages}{013133}.
%Type = Article
\bibitem[{Bush and Oza(2020)}]{Bush_2021}
\bibinfo{author}{J.~W.~M. Bush}, \bibinfo{author}{A.~U. Oza},
\newblock \bibinfo{title}{Hydrodynamic quantum analogs},
\newblock \bibinfo{journal}{Rep. Prog. Phy.} \bibinfo{volume}{84}
  (\bibinfo{year}{2020}) \bibinfo{pages}{017001}.
%Type = Article
\bibitem[{Bacot et~al.(2019)Bacot, Perrard, Labousse, Couder, and
  Fort}]{PhysRevLett.122.104303}
\bibinfo{author}{V.~Bacot}, \bibinfo{author}{S.~Perrard},
  \bibinfo{author}{M.~Labousse}, \bibinfo{author}{Y.~Couder},
  \bibinfo{author}{E.~Fort},
\newblock \bibinfo{title}{Multistable free states of an active particle from a
  coherent memory dynamics},
\newblock \bibinfo{journal}{Phys. Rev. Lett.} \bibinfo{volume}{122}
  (\bibinfo{year}{2019}) \bibinfo{pages}{104303}.
%Type = Article
\bibitem[{{Faraday}(1831)}]{Faraday1831a}
\bibinfo{author}{M.~{Faraday}},
\newblock \bibinfo{title}{On a peculiar class of acoustical figures; and on
  certain forms assumed by groups of particles upon vibrating elastic
  surfaces},
\newblock \bibinfo{journal}{Philosophical Transactions of the Royal Society
  London Series I} \bibinfo{volume}{121} (\bibinfo{year}{1831})
  \bibinfo{pages}{299--340}.
%Type = Article
\bibitem[{Durey et~al.(2020)Durey, Turton, and Bush}]{Durey2020}
\bibinfo{author}{M.~Durey}, \bibinfo{author}{S.~E. Turton},
  \bibinfo{author}{J.~W.~M. Bush},
\newblock \bibinfo{title}{Speed oscillations in classical pilot-wave dynamics},
\newblock \bibinfo{journal}{Proc. R. Soc. A} \bibinfo{volume}{476}
  (\bibinfo{year}{2020}) \bibinfo{pages}{20190884}.
%Type = Article
\bibitem[{Valani et~al.(2021)Valani, Slim, Paganin, Simula, and
  Vo}]{Valaniunsteady2021}
\bibinfo{author}{R.~N. Valani}, \bibinfo{author}{A.~C. Slim},
  \bibinfo{author}{D.~M. Paganin}, \bibinfo{author}{T.~P. Simula},
  \bibinfo{author}{T.~Vo},
\newblock \bibinfo{title}{Unsteady dynamics of a classical particle-wave
  entity},
\newblock \bibinfo{journal}{Phys. Rev. E} \bibinfo{volume}{104}
  (\bibinfo{year}{2021}) \bibinfo{pages}{015106}.
%Type = Phdthesis
\bibitem[{Mol{\'a}{\v{c}}ek(2013)}]{phdthesismolacek}
\bibinfo{author}{J.~Mol{\'a}{\v{c}}ek}, \bibinfo{title}{Bouncing and walking
  droplets : towards a hydrodynamic pilot-wave theory}, Ph.D. thesis,
  Massachusetts Institute of Technology, \bibinfo{year}{2013}.
%Type = Article
\bibitem[{Durey(2020)}]{Durey2020lorenz}
\bibinfo{author}{M.~Durey},
\newblock \bibinfo{title}{Bifurcations and chaos in a {L}orenz-like pilot-wave
  system},
\newblock \bibinfo{journal}{Chaos} \bibinfo{volume}{30} (\bibinfo{year}{2020})
  \bibinfo{pages}{103115}.
%Type = Article
\bibitem[{Valani(2022)}]{Valanilorenz2022}
\bibinfo{author}{R.~N. Valani},
\newblock \bibinfo{title}{{L}orenz-like systems emerging from an
  integro-differential trajectory equation of a one-dimensional wave–particle
  entity},
\newblock \bibinfo{journal}{Chaos} \bibinfo{volume}{32} (\bibinfo{year}{2022}).
%Type = Book
\bibitem[{Guenther and Lee(2012)}]{guentherlee2012}
\bibinfo{author}{R.~Guenther}, \bibinfo{author}{J.~Lee},
  \bibinfo{title}{Partial Differential Equations of Mathematical Physics and
  Integral Equations}, Dover Books on Mathematics, \bibinfo{publisher}{Dover
  Publications}, \bibinfo{year}{2012}.
%Type = Article
\bibitem[{Cox and Matthews(2002)}]{CM2002}
\bibinfo{author}{S.~Cox}, \bibinfo{author}{P.~Matthews},
\newblock \bibinfo{title}{Exponential time differencing for stiff systems},
\newblock \bibinfo{journal}{J. Comput. Phys.} \bibinfo{volume}{2}
  (\bibinfo{year}{2002}) \bibinfo{pages}{430--455}.
%Type = Article
\bibitem[{Mitchell and Vynnycky(2009)}]{mitchell2009}
\bibinfo{author}{S.~Mitchell}, \bibinfo{author}{M.~Vynnycky},
\newblock \bibinfo{title}{Finite-difference methods with increased accuracy and
  correct initialization for one-dimensional {S}tefan problems},
\newblock \bibinfo{journal}{Appl. Math. Comput.} \bibinfo{volume}{215}
  (\bibinfo{year}{2009}) \bibinfo{pages}{1609--1621}.

\end{thebibliography}
\end{document}